\documentclass[reqno,11pt]{amsart}
\usepackage{amsmath, latexsym, amsfonts, amssymb, amsthm, amscd}
\usepackage{graphics,epsf,psfrag}
\numberwithin{equation}{section}

\newtheorem{theorem}{Theorem}[section]
\newtheorem{lemma}[theorem]{Lemma}
\newtheorem{proposition}[theorem]{Proposition}

\newtheorem{rem}[theorem]{Remark}
\newtheorem{definition}[theorem]{Definition}

\DeclareMathOperator{\sign}{\mathrm{sign}}

\newcommand{\ind}{\mathbf{1}}
\newcommand{\dd}{\text{d}}
\newcommand{\R}{\mathbb{R}}

\newcommand{\N}{\mathbb{N}}
\renewcommand{\tilde}{\widetilde}

\newcommand{\cN}{{\ensuremath{\mathcal N}} }
\newcommand{\cL}{{\ensuremath{\mathcal L}} }

\newcommand{\cD}{{\ensuremath{\mathcal D}} }
\newcommand{\cS}{{\ensuremath{\mathcal S}} }
\newcommand{\cZ}{{\ensuremath{\mathcal Z}} }

\newcommand{\bP}{{\ensuremath{\mathbf P}} }
\newcommand{\bE}{{\ensuremath{\mathbf E}} }

\DeclareMathSymbol{\leqslant}{\mathalpha}{AMSa}{"36} 
\DeclareMathSymbol{\geqslant}{\mathalpha}{AMSa}{"3E} 
\DeclareMathSymbol{\eset}{\mathalpha}{AMSb}{"3F}     
\newcommand{\sumtwo}[2]{\sum_{\substack{#1 \\ #2}}} 
\newcommand{\prodtwo}[2]{\prod_{\substack{#1 \\ #2}}}     

\newcommand{\bbE}{{\ensuremath{\mathbb E}} }

\newcommand{\bbP}{{\ensuremath{\mathbb P}} }

\newcommand{\ga}{\alpha}
\newcommand{\gb}{\beta}
\newcommand{\gga}{\gamma}            
\newcommand{\gd}{\delta}
\newcommand{\gep}{\varepsilon}       
\newcommand{\gp}{\varphi}

\newcommand{\gD}{\Delta}

\newcommand{\go}{\omega}

\newcommand{\gl}{\lambda}

\newcommand{\gS}{\Sigma}

\makeatletter
\def\captionfont@{\footnotesize}
\def\captionheadfont@{\scshape}

\long\def\@makecaption#1#2{%
  \vspace{2mm}
  \setbox\@tempboxa\vbox{\color@setgroup
    \advance\hsize-6pc\noindent
    \captionfont@\captionheadfont@#1\@xp\@ifnotempty\@xp
        {\@cdr#2\@nil}{.\captionfont@\upshape\enspace#2}%
    \unskip\kern-6pc\par
    \global\setbox\@ne\lastbox\color@endgroup}%
  \ifhbox\@ne 
    \setbox\@ne\hbox{\unhbox\@ne\unskip\unskip\unpenalty\unkern}%
  \fi
  \ifdim\wd\@tempboxa=\z@ 
    \setbox\@ne\hbox to\columnwidth{\hss\kern-6pc\box\@ne\hss}%
  \else 
    \setbox\@ne\vbox{\unvbox\@tempboxa\parskip\z@skip
        \noindent\unhbox\@ne\advance\hsize-6pc\par}%
\fi
  \ifnum\@tempcnta<64 
    \addvspace\abovecaptionskip
    \moveright 3pc\box\@ne
  \else 
    \moveright 3pc\box\@ne
    \nobreak
    \vskip\belowcaptionskip
  \fi
\relax
}
\makeatother
\def\writefig#1 #2 #3 {\rlap{\kern #1 truecm
\raise #2 truecm \hbox{#3}}}

\newcommand{\tf}{\textsc{f}}
\newcommand{\M}{\textsc{M}}

\newcommand{\rc}{\mathtt c}

\begin{document}

\title[On the critical curve of copolymer models]{
Copolymers at 
selective interfaces:\\
new bounds on the phase diagram 
}

\author{Thierry Bodineau}
\address{D\'epartement de Math\'ematiques et Applications,
Ecole Normale Sup\'erieure,
 45 rue d'Ulm,
75230 PARIS Cedex 05, France
}
\email{bodineau\@@dma.ens.fr}

\author{Giambattista Giacomin}
\address{
  Universit{\'e} Paris 7 -- Denis Diderot and Laboratoire de Probabilit{\'e}s et Mod\`eles Al\'eatoires (CNRS U.M.R. 7599),
U.F.R.                Math\'ematiques, Case 7012,
                2 place Jussieu, 75251 Paris cedex 05, France
}
\email{giacomin\@@math.jussieu.fr}
\author{Hubert Lacoin}
\address{
  Universit{\'e} Paris 7 -- Denis Diderot and Laboratoire de Probabilit{\'e}s et Mod\`eles Al\'eatoires (CNRS U.M.R. 7599),
U.F.R.                Math\'ematiques, Case 7012,
                2 place Jussieu, 75251 Paris cedex 05, France
}
\email{lacoin\@@clipper.ens.fr}

\author{Fabio~Lucio~Toninelli}
\address{
Ecole Normale Sup\'erieure de Lyon, Laboratoire de Physique and CNRS,
UMR 5672, 46 All\'ee d'Italie, 69364 Lyon Cedex 07, France}
\email{fltonine@ens-lyon.fr}
\date{\today}

\begin{abstract} We investigate the phase diagram of disordered
  copolymers at the interface between two selective solvents, and in
  particular its weak-coupling behavior, encoded in the slope $m_c$ of
  the critical line at the origin. We focus on the directed walk case,
  which has turned out to be, in spite of the apparent simplicity,
  extremely challenging.  In mathematical terms, the partition
  function of such a model does not depend on all the details of the
  Markov chain that models the polymer, but only on the time elapsed
  between successive returns to zero and on whether the walk is in the
  upper or lower half plane between such returns. This observation
  leads to a natural generalization of the model, in terms of
  arbitrary laws of return times: the most interesting case being the
  one of return times with power law tails (with exponent $1+\ga$,
  $\ga=1/2$ in the case of the symmetric random walk).  The main
  results we present here are:
  \begin{enumerate}
  \item the improvement of the known  result $1/(1+\alpha)\le m_c\le 1$, as soon 
  as $\ga >1$ for what concerns the upper bound, and down to $\ga\approx 0.65$
  for the lower bound. 
  \item a proof of the fact that   the critical curve lies strictly below 
  the critical curve of the annealed model 
   for every non-zero value of the coupling parameter. 
  \end{enumerate}
  We also  
  provide an argument that rigorously shows the strong dependence
  of the phase diagram on the details of the return probability (and not
  only on the  tail behavior).
 Lower bounds are obtained by exhibiting a new localization strategy, while
   upper bounds are 
based on estimates  of non-integer moments of the partition function.
\\ \\ 2000 \textit{Mathematics
    Subject Classification:  60K35, 82B44, 60K37 } \\ \\
  \textit{Keywords: Directed Polymers, Disorder,  Copolymers at Selective Interfaces, Rare-stretch
    Strategies, Fractional Moment Estimates}
\end{abstract}

\maketitle

\section{Introduction }

Copolymers (or heteropolymers) are chains of non-identical monomer
units. We focus here on the case in which some of the monomer units
have an affinity for a solvent A, while the affinity of the others is
for a solvent B. Affinities - below we will call them {\sl charges} -
are fixed along the polymer chain and we will model them as quenched
disorder.  The medium in which the (co)polymer fluctuates is the one
schematized in Figure~\ref{fig:1}: the two solvents occupy half of the
space and they are separated by a sharp (and flat) interface.
Copolymer models have an extended literature, notably
models based on self-avoiding walks have been studied (see {\sl e.g.}
\cite{cf:Whittington} and references therein), but a very simple model,
that turned out to be nevertheless extremely challenging, has been
proposed in \cite{cf:GHLO}. It is a two-dimensional, in fact a
$(1+1)$-dimensional model in which the self-avoidance property is
enforced by considering directed walks and of course the walk steps
are the monomer units.  The $n$-th monomer of the chain carries a {\sl
  random charge} $\epsilon_n$ (which corresponds to $(\go_n+h)$ in
formula \eqref{eq:cZ} below), which can be either positive or
negative. Here {\sl random} refers to the fact that the charges are
not placed in a homogeneous or periodic way along the chain, but they
are (a realization of) a collection of independent and identically
distributed (IID) random variables.  The charge $\epsilon_n$
quantifies the chemical affinity of the $n$-th monomer with the
solvents so that the monomer has an energetic preference for being
placed in the solvent A (for example, oil) if $\epsilon_n>0$ and in
solvent B (water) if $\epsilon_n<0$.  It is then intuitive that an
energy-entropy competition takes place - maximizing the energy by
placing as many monomers as possible into their preferred solvent,
versus wandering away from the interface and gaining in entropy - and
this leads to a (non-trivial, as we will see)
localization-delocalization transition, when the temperature or the
mean of $\epsilon_n$ is varied.

The localization/delocalization critical curve $h_c(\gl)$, in the
$(\gl,h)$ plane (cf. Proposition \ref{th:hc}: $\gl$ is the coupling
parameter and $h$ is an asymmetry parameter controlling the mean of
the disorder), has attracted much attention, both in the theoretical
physics \cite{cf:CW,cf:GHLO,cf:Monthus,cf:SW,cf:SSE, cf:MT} and in
the mathematics literature
\cite{cf:AZ,cf:BisdH,cf:BG,cf:BdH,cf:Book,cf:GT,cf:Sinai}.  A point
that has to be stressed is that one can find in the physical
literature some predictions on the phase diagram and notably one can
find expressions, claimed to be exact, of the critical curve (see
below for more precision on this).  However, these expressions do not
coincide and this fact strongly suggests that the understanding of
this transition is very limited (this is confirmed by the lack of
agreement on critical exponents too
\cite{cf:CW,cf:Habi,cf:Monthus,cf:MT}).  In addition, the arguments in
many of the theoretical physics papers that we have mentioned appear
not to depend on the details of the return probabilities, but just on
the tail, in fact the arguments are developed for some large-scale
{\sl equivalent} systems in which the local details are forgotten. One
of the results we present below shows that the critical curve does
depend on these details, and at times even in a very radical way.
What (most probably) does not depend on the details of the polymer and
of the disorder is the slope $m_c$ of the critical curve at the
origin, which corresponds to weak polymer-solvent coupling.

Therefore $m_c$ is definitely a quantity of great interest for its
{\sl universal} character. 

From a mathematical standpoint  bounds on the
critical line have been established in 
 \cite{cf:BG,cf:BdH}, but they are not sufficiently precise to 
 settle the controversies between  the different
physical predictions (a result on the critical behavior has been
proven in \cite{cf:GT_cmp}; it shows that the transition is in {\sl
  great generality} at least of second order). These bounds in some
very exceptional cases do coincide (but these are really marginal
situations, not considered by physicists). 
In particular the upper
bound follows  from  an annealing procedure 
 \cite{cf:BdH}.
Very recently one of us \cite{cf:T_fractmom} has proven that at
large coupling parameter and for unbounded disorder the  bound
in \cite{cf:BdH}
can be improved. Such a result is based on estimates on the fractional
moments of the partition function.  Here we will go beyond this result
and get upper bounds that hold for arbitrary coupling parameter and 
general charge distributions, still
by estimating fractional moments and by adapting an idea first
developed in the context of disordered pinning models \cite{cf:DGLT}.
Let us mention that a result close to ours
 has been recently obtained, using  different
techniques,  by E. Bolthausen and
F. den Hollander \cite{cf:BdH2}.
It should be pointed out that, for copolymer models, constrained
annealing techniques have been applied at several instances
({\sl e.g.} \cite{cf:IRW} and references therein), but it has been shown in \cite{cf:CG} that
they are useless to go beyond the  bound in \cite{cf:BdH} on the
critical line.

Our purpose is to improve also the lower bounds.  Finding lower bounds
on $h_c(\gl)$ amounts to finding lower bounds on the partition
function, and to this aim it is natural to try to guess what the most
favorable polymer configurations are and to keep only those in the
partition sum. We will refer to such a trajectory selection as to a
(selection) {\sl strategy}.  However natural this idea may look, it is
difficult to guess strategies which give non-trivial bounds.  In
\cite{cf:Monthus}, a real-space renormalization procedure was
implemented by using a {\sl rare-stretch strategy} which takes
advantage of atypical regions in the disorder.  For some time we have
believed that the strategy of \cite{cf:Monthus} would yield a correct
description of the critical curve and the correct value of $m_c$, even
if in mathematical terms such a strategy led only to a lower bound
\cite{cf:BG}.  Such a belief was later shaken by accurate numerical
simulations \cite{cf:CGG}.  In this paper we present a new
rare-stretch strategy which, although not optimal, improves in some
situations the critical curve lower bounds \cite{cf:BG} based on the
renormalization approach of \cite{cf:Monthus}.  \smallskip

\begin{rem}
  \rm In the physical literature only models based on symmetric walks
  have been considered. However we find that considering more general
  models helps in a substantial way in devising new arguments of proof
  and in understanding the limitations of previous approaches.
  Besides, the {\sl generalized copolymer model} (already introduced
  and studied in \cite{cf:GTloc,cf:Book}) is a natural, and very
  easily defined, disordered model where one can investigate the
  effect of disorder on systems for which the  associated annealed model
  has a first order phase transition.
\end{rem}

Finally, we mention that also the situation where the geometry of the
regions occupied by the two solvents is more involved than just two
half-planes has been considered, see for instance \cite{cf:PdH}.

\section{The model and the results}

\subsection{The standard copolymer model}
\label{sec:standard}
\begin{figure}[h]
\begin{center}
\leavevmode
\epsfysize =6.8 cm
\psfragscanon
\psfrag{0}[c]{$0$}
\psfrag{n}[c]{$n$}
\psfrag{Sn}[c]{$S_n$}
\psfrag{sA}[c]{\small Solvent A}
\psfrag{sB}[c]{\small Solvent B}
\psfrag{Int}[c]{\small Interface}
\psfrag{o1}[c]{\small $\go_1$}
\psfrag{o2}[c]{\small $\go_2$}
\psfrag{o3}[c]{\small $\go_3$}
\psfrag{o4}[c]{\small $\go_4$}
\psfrag{t1}[c]{$\tau_1$}
\psfrag{t2}[c]{$\tau_2$}
\psfrag{t3}[c]{$\tau_3$}
\psfrag{t4}[c]{$\tau_4$}
\epsfbox{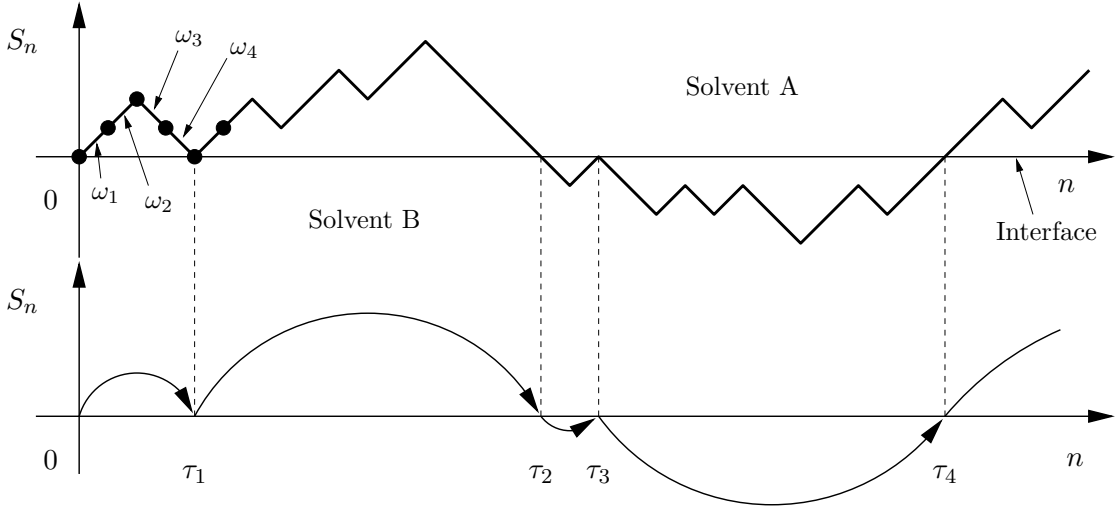}
\end{center}
\caption{\label{fig:1} In the top drawing  one finds a trajectory of the
  standard copolymer model. The simple random walk bonds, that is the
  segments linking $(n-1, S_{n-1})$ and $(n, S_{n})$, are the monomers
  and they carry a charge $\go_n$ that is drawn at random, but it is
  fixed once for all ({\sl quenched} disorder), while the polymer
  fluctuates. Positively charged monomers are energetically rewarded
  if they lie in the upper half-plane, occupied by solvent A, while
  they are penalized if they lie in the lower half-plane (solvent B).
  The situation is reversed for negatively charged monomers. If the
  parameter $h$ is not zero, the model is not symmetric under the
  exchange of the solvents. Note that the energetically favored
  trajectories are the ones that place {\sl most of} the monomers in
  their favored solvent, but such trajectories are necessarily
  sticking close to the interface between the two solvents and they
  are therefore few compared to the trajectories that wander more
  freely. This leads to an energy-entropy
  competition between localized and delocalized trajectories.  Observe
  also that the energy does not depend on the details of the
  trajectory between successive visits to $0$, so that the model
  can be schematized like in the lower figure, that is simply in terms
  of $\tau_1=2$, $\tau_2-\tau_1=6$, $\tau_3-\tau_2=1$,
  $\tau_4-\tau_3=6$,... (that are distances between successive returns to
  zero of $S_{2n}$) and of the sign sequence $s_1=+1$, $s_2=+1$,
  $s_3=-1$, $s_4=-1$, $s_5=+1$,... The schematized version of the
  model can be easily generalized to arbitrary return laws $K(\cdot)$.  }
\end{figure}

The model of copolymers close to a selective interface introduced in
\cite{cf:GHLO} is based on the simple random walk $S:=\{S_n\}_{n=0,1,
  \ldots}$, that is the Markov chain characterized by $S_0=0$ and the
fact that the increment sequence $\{S_{n+1}-S_n\}_{n=0,1, \ldots}$ is
IID with $\bP(S_1=+1)= \bP(S_1=-1)=1/2$. The partition function of the
model is defined as
\begin{equation}
\label{eq:cZ}
\cZ_{N, \go} \, := \, 
\bE \left[ \exp \left(\gl \sum_{n=1}^N (\go_n +h) \sign (S_n)
\right)
\right],
\end{equation} 
where $N$ is a positive integer, $\bE$ is the expectation with respect 
to the random walk trajectory, $\gl$ and $h$ are two constants that can be chosen non-negative without loss of generality, and $\go := \{ \go_n \}_{n=1,2, \ldots}$ is a sequence 
of real numbers. Since $\sign(0)$ is {\sl a priori}  not defined we stipulate
that $\sign (S_n)= \sign(S_{n-1})$ if $S_n=0$ (this arbitrary choice appears  natural 
once the process is {\sl decomposed into excursions}, see  below).
Let us also remark that the symbol $\bE$ just denotes the normalized sum
over the $2^N$ trajectories of the simple random walk. Figure 
\ref{fig:1} may be of help
in order to get some insight.

\smallskip

\begin{definition} \label{def:charges}\rm
Unless otherwise stated,
the sequence $\go$, referred to as {\sl sequence of charges}, is chosen as a typical
realization of an IID sequence of law $\bbP$. We assume  $ \M (t):=\bbE[ \exp(t \go_1)]< \infty$
for every $t \in \R$,  and that $\bbE[\go_1]=0$ and $\bbE [\go_1^2]=1$.
\end{definition}
\smallskip

A useful observation about this model is that $\cZ_{N, \go}$ can be
expressed in terms of the return times
$\tau:=\{\tau_j\}_{j=0,1,\ldots}$ defined iteratively by $\tau_0=0$
and $\tau_{j+1}:= \inf\{ n> \tau_j: \, S_n=0\}$.  Note that $\tau$ is
a random walk itself: it has (positive and) IID increments $\{
\tau_{j+1}-\tau_j\}_{j=0,1, \ldots}$ with a law $K(n)=\bP(\tau_1=n)$
that is explicitly written in combinatorial terms (and, notably,
$\lim_{n \to \infty} n^{3/2}K(2n)=\sqrt{1/(4\pi)}$). By using a
standard (probabilistic) terminology, we say that $\tau$ is a renewal
sequence.  Since $\sign(S_n)$ is constant inside an {\sl excursion}
$\{\tau_j +1, \tau_j +2, \ldots, \tau_{j+1}\}$, it is natural to
consider the sequence defined by $s_j = \sign \left( S_{\tau_j
  }\right), j\ge1$.  An immediate consequence of the (strong) Markov
property is that $s= \{s_j\}_{j=1,2, \ldots}$ is an IID sequence of
symmetric random variables (taking of course only the values $\pm 1$).

Before generalizing the model, let us immediately say that it is technically
advantageous to work with a slightly different definition of the
energy (and of the partition function \eqref{eq:cZ}) of the model given by
\begin{equation}
\label{eq:Z}
Z_{N, \go} \, := \, 
\bE \left[ \exp \left(-2\gl \sum_{n=1}^N (\go_n +h) \gD_n
\right)
\right],
\end{equation}
where $\gD_n:= [1-\sign(S_n)]/2$, so $\gD_n$ is the indicator function
that the $n^{\rm th}$-monomer is below the interface.
Since the difference between the terms in the 
exponent of the expressions in \eqref{eq:cZ} and \eqref{eq:Z}
is just $\gl \sum_{n=1}^N (\go_n +h)$ (in particular, independent of $S$)
the two models are actually the same and the asymptotic behaviors of 
$\cZ_{N, \go}$ and $Z_{N, \go}$ are trivially related.

\subsection{The generalized copolymer model}
\label{sec:general}

The observation we have just made naturally leads to a generalization
of the framework.  
We are in fact going to assume that $K(\cdot)$ is a
discrete probability density on $\N$, so that $\sum_{n\in\N}
K(n)=1$, such that for some $\ga>0$
\begin{equation}
\label{eq:K}
K(n)\stackrel{n \to \infty}\sim \frac{C_K}{n^{1+\ga}},
\end{equation}
where $a_n \stackrel{n \to \infty}\sim b_n$ means that $\lim_{n \to
  \infty} a_n/b_n =1$, and $C_K$ is a positive constant.  We define
then a renewal process $\tau:=\{\tau_0, \tau_1, \ldots\}$, {\sl i.e.} a random walk with positive
(and independent) increments such that $\tau_0=0$ and the sequence $\{
\tau_{j}-\tau_{j-1}\}_{j=1, 2, \ldots}$ is an IID sequence with law
$\bP (\tau_1=n)=K(n)$ for every $n$.  We also let $s= \{s_j\}_{j}$ be
an IID sequence of symmetric random variables taking values $\pm 1$,
and for $n\in\N$ we define $\Delta_n=(1-s_j)/2$ if $\tau_{j-1}< n\le \tau_j$.
The partition function of the generalized model is then again \eqref{eq:Z},
where now $\bE$ denotes the expectation with respect to $\tau$ and $s$.

Note that, strictly speaking, \eqref{eq:K} is not compatible with the
random walk choice which has only even return times. However, one can
always redefine the return times of the simple random walk as
$\tau/2$, and this is of course a trivial change (and $\ga=1/2$).

\medskip

\begin{rem}
\label{rem:L}\rm 
There is no difficulty in relaxing the assumption \eqref{eq:K}, for example
by allowing {\sl logarithmic} corrections to the asymptotic behavior. 
For the sake of simplicity, we are going to stick to assumption \eqref{eq:K},
with the notable exception of Remark~\ref{rem:1} and \S~\ref{sec:further}.
\end{rem}

\subsection{The free energy and the phase diagram}
\label{sec:fe}

The free energy (per unit length) of such a model is
\begin{equation}
\label{eq:fe}
\tf (\gl, h) \, :=\, \lim_{N \to \infty} \frac 1N \bbE \log Z_{N , \go}.
\end{equation}
Actually, the $\bbP(\dd \go)$-a.s. limit of the sequence
of random variables $\{  N^{-1}\log Z_{N , \go}\}_N$ exists
and coincides with $\tf (\gl, h)$ (see, {\sl e.g.}, \cite[Ch.~4]{cf:Book}). 
We observe that $\tf (\gl, h)\ge 0$ because
\begin{multline}
\label{eq:lb0}
 Z_{N , \go} \, \ge \,  
 \bE \left[ \exp \left(-2\gl \sum_{n=1}^N (\go_n +h) \gD_n
\right); \, \gD_n =0 \text{ for } n=1,2, \ldots, N
\right]
\\
=  
 \bP \left(  \gD_n =0, \,  n=1,2, \ldots, N
\right)\, =\, \bP( s_1=+1, \tau_1\ge N) \, =\, \frac 12
\sum_{n \ge N } K(n)\stackrel{N \to \infty}\sim
\frac{C_K}{2\ga  N^{\ga}}.
\end{multline}
In more intuitive terms we say that the free energy in the delocalized regime
is zero and we split the phase diagram according to 
\begin{equation}
\cL\,=\, \left\{ (\gl, h): \, \tf (\gl, h)>0\right\} \ \ \text{ and  } \ Ê\ 
\cD\,=\, \left\{ (\gl, h): \, \tf (\gl, h)=0\right\}.
\end{equation}

\medskip

\begin{rem}
\label{rem:superadd} \rm 
Not surprisingly, the proof of  the existence of the limit in
\eqref{eq:fe} relies on super-additivity. Super-additivity turns out to be a very 
crucial tool for our arguments too, so let us stress that 
$\left\{\bbE \log Z_{N, \go}\right\}_N$ is {\sl not} a super-additive sequence. Rather one has to consider
\begin{equation}
\label{eq:Zc}
Z_{N, \go}^{\rc} \, := \, 
\bE \left[ \exp \left(-2\gl \sum_{n=1}^N (\go_n +h) \gD_n
\right); \, N\in\tau
\right],
\end{equation}
(${\rc}$ for {\sl constrained})
where  $N\in\tau$ simply means that $\tau_n=N$ for some $n$.  It
is easy to see that $\left\{\bbE \log Z_{N, \go}^\rc\right\}_N$ is
super-additive (and it is just a bit harder to see that $ (1/N)\bbE
\log Z_{N, \go}^\rc -(1/N) \bbE \log Z_{N, \go}\to 0$ as $N \to \infty$, so
that, when talking of the free energy, we can safely switch between
the two partition functions).  Super-additivity says also that $\tf
(\gl, h)$ coincides with $\sup_N N^{-1 } \bbE \log Z_{N, \go}^\rc$ and
this yields the following characterization:
\begin{equation}
\label{eq:criterion}
\tf (\gl,h)>0  \text{ if and only if there exists } N \text{ such that }
 \bbE \log Z_{N, \go}^\rc>0.
\end{equation}
This is a powerful tool because it is a finite-volume criterion for
localization. It has played a central role in a number of results,
like \cite{cf:CGG,cf:GT}, and it will be, again, very important here. 
\end{rem}

\medskip

\begin{rem}
  \rm In this work we leave aside any consideration on the copolymer
  path behavior and concentrate on free-energy properties.  The fact
  that $(\gl, h)\in \cL $ (respectively $(\gl, h)\in \cD $) does
  correspond to localized (respectively, delocalized) behavior of the
  paths of the process has been addressed in depth elsewhere ({\sl
    e.g.}, \cite{cf:BisdH,cf:GTloc,cf:GT} and \cite[Ch.s~7 and 8]{cf:Book}).
\end{rem}
\medskip

Convexity and monotonicity properties of the free energy entail
a number of properties of the phase diagram that we sum up
in the next statement (see \cite[Ch.~6]{cf:Book} for a proof, built on
results proven in \cite{cf:BdH,cf:BG}).

\medskip

\begin{proposition}
\label{th:hc} ({\sl Existence of the critical curve}.)
There exists a continuous, strictly increasing function
$\gl\mapsto h_c(\gl)$, satisfying $h_c(0)=0$, such that
$\cD= \{(\gl, h):\, h \ge h_c(\gl)\}$ (and of course $\cL =\{(\gl, h):\, h <h_c(\gl)\}$).
\end{proposition}

\medskip

This note addresses precisely the behavior of $h_c(\cdot)$.
 Let us first recall the known results: in the next statement 
we give bounds  
 proven for the standard copolymer model
in \cite{cf:BdH} for what concerns the upper bound and in \cite{cf:BG}
for the lower bound (the straightforward adaptation of the arguments
to cover the generalized model is detailed in \cite[Ch.~6]{cf:Book}).
We set for $m>0$
\begin{equation}
\label{eq:hm}
h^{(m)}(\gl)\,:=\, \frac 1{2m\gl} \log \M (-2m\gl ),
\end{equation}
where
$\M(t)$ is given in Definition~\ref{def:charges} (note that $\dd h^{(m)}(\gl)/\dd \gl \vert _{\gl=0}=m$).

\medskip

\begin{proposition}
\label{th:ublb}
For every choice of $K(\cdot)$ satisfying \eqref{eq:K} we have
\begin{equation}
\label{eq:ublb}
h^{(1/(1+\ga))}(\gl) \, \le \, h_c(\gl) \, \le \,h^{(1)}(\gl),
\end{equation}
for every $\gl \ge 0$.
As a consequence,
\begin{equation}
\label{eq:slopelb}
\frac 1{1+\ga} \, \le \, 
\liminf_{\gl \searrow 0} \frac{h_c(\gl)} \gl \, \le \, 
\limsup_{\gl \searrow 0} \frac{h_c(\gl)} \gl \, \le \, 1.
\end{equation}
\end{proposition}
\medskip

\begin{rem}\label{rem:1}\rm

Upper and lower bounds in \eqref{eq:ublb}
coincide only if $\ga=0$, but assumption \eqref{eq:K} requires $\ga >0$ 
in order to ensure that $K(\cdot)$ is normalizable.
 We can set $\ga=0$ if we accept to relax 
somewhat \eqref{eq:K} ({\sl cf.}
 Remark~\ref{rem:L}) and if we choose for example
 $K(n)=c/(n (\log n)^2)$ ($c$ such that $\sum_n K(n)=1$). 
 Since \eqref{eq:ublb}  holds also when there are logarithmic corrections
 to the power law behavior of $K(\cdot)$, and in particular for $\ga=0$
  (cf. \cite[Ch.~6]{cf:Book}), it is straightforward to see that $h_c(\gl)= h^{(1)}(\gl)$
 for every $\gl \ge 0$.
 \end{rem}

\subsection{Weak-coupling limit, rare-stretch strategy and a look at the literature} 
\label{sec:rev}

Much work has been done on the 
copolymer model, both in the physical and mathematical literature.
In spite of this, the understanding of the model is still very limited.
An important point  has been set forth in \cite{cf:BdH}, where
it has been shown, for the case of the standard model 
(in particular, $\ga=1/2$) and for $\go_1$ binary random variable,
that 
\begin{equation}
  \label{eq:convF}
  \lim_{\gga \searrow 0} \frac 1{\gga^2}
  \tf (\gga\gl ,\gga h)\, =\,  \phi (\gl, h),
\end{equation}
where
\begin{equation}
\label{eq:phi}
\phi(\gl,h) \,:=\, \lim_{t \to \infty} \frac 1t \bbE \log 
\bE \left[ \exp \left( 
-2\gl \int_0^t \ind_{B(t) <0} \left( \dd \gb (t) +h \dd t  \right) \right)
\right],
\end{equation} 
and $B$ and $\gb$ are two standard Brownian motions, respectively of
law $\bP$ and $\bbP$. Note that \eqref{eq:phi} is the partition
function of a continuous copolymer model; like for the discrete case
it is easy to see that $\phi(\gl, h)\ge 0$ and once again one can
define a localized and a delocalized regime, which are separated by a
continuous critical curve $\gl\mapsto \tilde h_c(\gl)$, {\sl i.e.}
the analog of Proposition~\ref{th:hc} holds. The novelty is that from the scaling properties of Brownian motion 
there exists a non-negative number $m_c$ such that $\tilde h_c (\gl)=m_c \gl$ for every $\gl\ge 0$.  
This is not all: in \cite{cf:BdH} it is shown that
\begin{equation}
\label{eq:slope}
\lim_{\gl \searrow 0} \frac{h_c(\gl)}\gl \, =\, m_c,
\end{equation}
a result which is more subtle than  the convergence
\eqref{eq:convF} of the free energy. 

\medskip

The results we have just stated, as well as their proofs, have a
strong flavor of {\sl universality}, in the sense that they are based
on the idea that at small coupling ($\gl $ small) typical excursions
are very long (since the simple random walk is null recurrent), and
small excursions do not contribute to the energy, so that the walk
can be replaced by a Brownian motion and the sum of the charges within
an excursion can be approximated by a normal variable. This does not
seem to be specific of the simple random walk and binary charges and
in fact the results have been shown to hold for much more general
charges \cite{cf:GT}. However, the results are (by far) not a direct
consequence of standard invariance principles and the
generalization to more general walks is highly non trivial.
 
It would be very interesting to extend the weak-coupling results
stated previously by simply assuming the validity of \eqref{eq:K} with
$\ga=1/2$. We expect in particular 
\eqref{eq:convF}-\eqref{eq:slope} to hold in such a generality. More
generally, one would like to have the weak-coupling results for
general $\ga \in (0,1)$.  Of course the expression \eqref{eq:phi} has
to be suitably changed, but we still expect \eqref{eq:slope} to hold
with some $m_c$ depending on $K(\cdot)$ only through $\alpha$.

In spite of the fact that the weak-coupling results we have stated do
not go yet as far as we would like, they clearly point the attention
to the slope of the critical curve at the origin as a quantity of
great interest.  And, at least for the case $\ga=1/2$, this issue has
been addressed in the physical literature, but without a consensus. In
particular in \cite{cf:GHLO} and \cite{cf:MT} it is claimed that
$m_c=1$ (note that $m_c\le 1$ by the upper bound in
\eqref{eq:slopelb}) while in \cite{cf:Monthus} and \cite{cf:SSE} it is
claimed that $m_c=2/3$ (and we know that $m_c\ge 2/3$ by the lower
bound in \eqref{eq:slopelb}).  For the standard copolymer model there
exists numerical evidence that $2/3<m_c < 1$, and possibly that
$m_c\approx 0.83$ \cite{cf:CGG} (see also \cite{cf:MG}), but until now
there is not much clue on how to estimate this value beyond
\eqref{eq:slopelb}.  Some of the papers we just mentioned are actually
claiming that $h_c(\gl)=h^{(m)}(\gl)$ for every $\gl$, with $m$ either
equal to $2/3$ or $1$.  However only \cite{cf:Monthus} deals with
general disorder laws, while the others are restricted to Gaussian
disorder, for which, incidentally, $h^{(m)}(\gl)=m\gl$.

\medskip

\begin{rem}
\label{rem:improve}
\rm In \cite{cf:T_fractmom} it has been shown that the upper bound in
\eqref{eq:ublb} can be improved for all $\ga$ if $\gl$ is
large and $\bbP(\go_1<x)>0$ for every $x$ (which is true for instance
in the case of Gaussian disorder). This is discussed briefly in
Section \ref{sec:upper} below. On the other hand, in
Proposition~\ref{th:further} below we show also that, again for all
$\ga$, we can find suitable inter-arrival distributions $K(\cdot)$ for
which the lower bound in Proposition~\ref{th:ublb} is not optimal.
These results however do not give any information on the slope of
$h_c(\cdot)$ at the origin.

\end{rem}  

\medskip

\begin{figure}[h!]
\begin{center}
\leavevmode
\epsfxsize =14 cm
\psfragscanon
\psfrag{0}[c]{$0$}
\psfrag{n}[c]{$n$}
\psfrag{Sn}[c]{$S_n$}
\psfrag{l}[c]{\small $\ell$}
\psfrag{b}[c]{\small $2\ell$}
\psfrag{c}[c]{\small $3\ell$}
\psfrag{d}[c]{\small $4\ell$}
\psfrag{e}[c]{\small $5\ell$}
\psfrag{f}[c]{\small $6\ell$}
\psfrag{g}[c]{\small $12\ell$}
\psfrag{h}[c]{\small $13\ell$}
\psfrag{i}[c]{\small $16\ell$}
\psfrag{m}[c]{\small $17\ell$}
\epsfbox{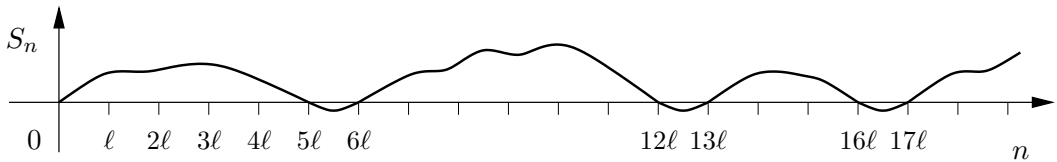}
\end{center}
\caption{\label{fig:2} The lower bound in \eqref{eq:ublb} is proven by
  using a coarse-graining parameter $\ell$ ($1 \ll \ell \ll N$). We
  consider the block charges $Q_j= \sum_{n=j\ell+1}^{(j+1)\ell}
  \go_n$, $j=0,1, 2,\ldots$, and we single out the blocks ({\sl rich
    blocks}) for which $ Q_j \le -q \ell$ ($q$ is a parameter that is
  going to be optimized).  The sequence $\{Q_j\}_j$ is of course IID,
  so that the location of the rich blocks is simply given by a
  Bernoulli trial sequence with parameter $p(\ell):= \bbP( Q_1 \le -q
  \ell)$.  For $\ell$ large, $p(\ell)$ is exponentially small, so that the
  rich blocks are rare and typically very spaced. 
  The lower bound in \eqref{eq:ublb} is obtained by bounding from below
  the free energy. This is achieved by restricting the partition function to the 
  trajectories $S$ visiting the
  lower half-plane only in correspondence of the rich blocks, {\sl
    i.e.}, such that $\sign(S_n)<0$ if and only if $n$ belongs to a
  rich block. 
These trajectories have very long excursions before returning to rich blocks as depicted 
in the drawing (the rich blocks are in this case
  the $5^{\rm th}$, the $12^{\rm th}$ and the $16^{\rm th}$). We refer
  to \cite{cf:BG,cf:Book} for a full proof of the lower bound in
  \eqref{eq:ublb}.  }
\end{figure}

For the sequel of the paper it is also important to sketch the idea
that leads to the lower bound in \eqref{eq:ublb}.  We do this with
Figure~\ref{fig:2} and its caption. In this strategy, the partition
function is evaluated only on trajectories which are made of very long
excursions with returns in {\sl rare stretches} where the mean of the
charges is atypical. This reduces the complexity of the model in two
ways. On  one hand, the charges act only as very rare  {\sl energetic
traps} and, on the other hand, the trajectories gain only the averaged
charge of these traps. A toy version of the copolymer mimicking this
behavior has been introduced in \cite{cf:BG} and it has been recently
proved \cite{cf:T_fractmom, cf:BCT} that in the strong-disorder limit
the simple rare-stretch strategy of Figure~\ref{fig:2} identifies
correctly the asymptotic behavior of its critical line.  However, the
numerical simulation we mentioned above show that this strategy
does not fully catch the correct behavior of the copolymer model
\eqref{eq:Z}.

\subsection{New results}

\subsubsection{Lower bound: a new rare-stretch strategy}
\label{sec:neutral}

A limitation in the rare-stretch strategy
described in Figure \ref{fig:2} is in the triviality of the behavior
of the polymer in the rare advantageous stretches (the {\sl rich
blocks}). Designing a better strategy is however not so obvious.
Here we are going to present one that is ultimately going back to the
original idea in \cite{cf:Sinai} that a neutral environment suffices
to localize the polymer (we will have to make this quantitative).  The
idea is therefore to look for neutral stretches, {\sl i.e.} $q=0$ (instead of negative stretches $-q<0$), and
employ a non-trivial localization strategy in these regions. The most
interesting results we have been able to extract from such idea are
summarized in the next statement.

\medskip

\begin{theorem}
\label{th:neutral}
For $\ga\ge 1$
\begin{equation}
\label{eq:neutral1}
\liminf_{\gl \searrow 0} \frac{h_c(\gl)} \gl \, \ge \, \max\left(\frac 1{\sqrt{1+\ga}},\frac12\right).
\end{equation}
Moreover $\ga \ge 0.801$ we have
\begin{equation}
\label{eq:neutral2}
\liminf_{\gl \searrow 0} \frac{h_c(\gl)} \gl \, > \, \frac 1{{1+\ga}}. 
\end{equation}
\end{theorem}

\medskip

Of course these results acquire a particular interest when compared
with \eqref{eq:slopelb}. We want to stress that the lower bound
on the slope that we are able to establish is rather explicit for all values of $\ga$
(see Section~\ref{sec:LB}, Proposition~\ref{th:frombook}),  in particular  
we are able to establish \eqref{eq:neutral2} for every
$\ga > \ga_0$, with 
$2(1+\ga_0)A(\ga_0)=1$, where $A(\ga)$ is the expression in the right-hand side of \eqref{eq:quasiexpl} optimized with respect to  $\kappa$. 
The explicit bound given in the statement comes from 
an explicit lower bound on  $A$, but
$A$ can be computed via numerical integration and one sees that $\ga_0$ is smaller than (but close to) $0.65$.
 In Remark~\ref{rem:nienteunmezzo}  we are going to explain that there is very little hope to make this strategy 
 work for $\ga=1/2$. Nevertheless, this shows that there are better strategies
 than the rare stretch strategy leading to the $1/(1+\ga)$ bound on the slope.

\subsubsection{Upper bound: fractional moments and the slope}
\label{sec:upper}
The upper bound in \eqref{eq:ublb} is just a consequence of the annealed inequality for the free energy:
\begin{eqnarray}
  \label{eq:annineq}
 \frac1N \bbE \log Z_{N,\go}  \le\frac1N \log \bbE Z_{N,\go}.
\end{eqnarray}
A natural idea to go beyond simple-minded annealing is to observe that,
again thanks to Jensen's inequality, for every $\gamma>0$
\begin{eqnarray}
  \label{eq:eghein}
 \frac1N \bbE \log Z_{N,\go}  \le\frac1{N\gamma}\log \bbE
\left[ (Z_{N,\go})^\gamma\right].
\end{eqnarray}
Since for $\gamma=1$ one recovers annealing \eqref{eq:annineq} and for
$\gamma\searrow 0$ \eqref{eq:eghein} becomes an equality, it is
natural to hope that non-trivial information can be obtained
estimating the $\gamma$-moments of the partition function, with
$0<\gamma<1$. This is precisely the approach which was followed in
\cite{cf:T_fractmom}, and indeed it turns out that at least for $\gl$
large enough and assuming that the random variables $\go_n$ are
unbounded (more precisely, that $\bbP(\go_1<x)>0$ for every $x$) one
can prove that $h_c(\gl)<h^{(1)}(\gl)$ (of course, this says nothing
about the critical slope). Let us also recall that in \cite[Corollary
3.9]{cf:T_fractmom} the same method allowed to prove that if
$\sum_{n\ge1}K(n)<1$ ($\tau$ is transient, a case that we are not
considering here) then $h_c(\gl)<h^{(1)}(\gl)$
{\sl for every} $\lambda>0$ and $\limsup_{\gl\searrow0}h_c(\gl)/\gl<1$.

Two important ingredients were added in Refs.  \cite{cf:GLT,cf:DGLT}
where, in the somewhat different context of {\sl disordered
  pinning/wetting models}, it was realized first of all that it is
actually sufficient to control the fractional moments of $Z_{N, \go}$
up to $N$ of the order of the correlation length of the annealed
system, and secondly that this control can be obtained through a
change-of-measure argument.  Here we generalize these arguments to the
case of the copolymer, and our main result is the following:

\medskip

\begin{theorem}
\label{th:ub}
Choose $K(\cdot)$ satisfying \eqref{eq:K}. 

If $\ga>1$ 
\begin{equation}
\label{eq:ub}
\limsup_{\gl \searrow 0} \frac{h_c(\gl)} \gl \,< \,1.
\end{equation}

If $0<\ga\le1$ there exists a positive constant $c$ such
that, for $0<\gl\le 1$
\begin{equation}
  \label{eq:a<1}
  h_c(\gl)\,\le\,  h^{(1)}\left(\gl \left(1-\frac c{\left|\log c\gl ^2\right|}\right)\right).
\end{equation}

Moreover, for every $\ga>0$ and $\gl>0$ one has $h_c(\gl)<
h^{(1)}(\gl)$.
\end{theorem}

\medskip

\subsubsection{A further remark on the critical curve}
\label{sec:further}

We complete the list of new results with the following one that
becomes of interest in view of the various conjectures that one finds
in the literature: in short, it says that in general $h_c(\cdot)$
heavily depends on the details of $K(\cdot)$ and it is certainly not
simply a function of $\ga$.

\medskip
\begin{proposition}
\label{th:further}
For every $\ga >0$, every $\gl>0$ and every $\gep>0$ we can find
$K(\cdot)$ that satisfies \eqref{eq:K} such that $h_c(\gl)>
h^{(1)}(\gl) -\gep$.
\end{proposition}
\medskip

The proof  is short and it is of help in understanding the result, so we give it right here. 
\smallskip

\noindent
{\it Proof.}
Consider the model with $\ga=0$ ({\sl cf.} Remark~\ref{rem:1}, call
$\tilde K(\cdot)$ the particular return probability of that model).
Fix $\gl>0$ and $\gep>0$: we have $\tf (\gl, h^{(1)}(\gl)-\gep)>0$,
so that, by \eqref{eq:criterion}, there exists $N\in \N$
such that $\bbE \log Z^\rc_{N, \go} >0$. 
But such a result is unchanged if we modify the definition
of $\tilde K(n)$ for $n >N$, since $Z^\rc_{N, \go} $ itself is unchanged. 
Since \eqref{eq:K} depends only on the tail of $K(\cdot)$, we are done.
\qed

\medskip

In its simplicity, Proposition~\ref{th:further}, in combination with
the other results we have stated or that can be found in the
literature, can be used to rule out a number of conjectures that have
been made or that one may be tempted to make. For example, when the
disorder is Gaussian the critical line in general is not a straight
line with a slope depending only on $\ga$
(recall Remark~\ref{rem:improve}).  More generally, it casts serious
doubts on the fact that the critical line coincides with
$h^{(m_c)}(\gl)$ for every $\gl$ regardless of the details of
$K(\cdot)$, even if numerical evidence in \cite{cf:CGG} suggested that
the critical line could coincide with $h^{(m_c)}(\cdot)$ in the
particular case of the standard copolymer model.

\medskip

For the remainder, we stress  that
we systematically develop the arguments first for $\go_1\sim \cN (0,1)$ and 
then give
 the modifications needed  to deal with the general charge
distributions of Definition~\ref{def:charges}.

\section{Neutral stretches and a lower bound strategy}
\label{sec:LB}

In this section we are going to give a proof of  Theorem~\ref{th:neutral}.

\begin{proposition}
\label{th:0strategy-Gauss}
Let us consider the general copolymer model with $\go_1\sim \cN (0,1)$.
For every $\gl>0$ we have 
\begin{equation}
\label{eq:0strategy-Gauss}
h_c(\gl) \, \ge \, \sqrt{\frac {2 \: \tf (\gl,0)} {1+\ga}}.
\end{equation} 
\end{proposition}
\medskip

\noindent
{\it Proof.}
For  $\ell\in \N$ let us define the random variable 
\begin{equation}
  \tf_\ell (\gl, h; \go)\, :=\, 
  \frac 1\ell \log \bE 
\left[ \exp\left(-2\gl \sum_{n=1}^\ell (\go+h) \gD_n\right) ; 
\ell \in \tau \right].
\end{equation}
We now claim that for every $\gd >0$ we have
\begin{equation}
\label{eq:claim-Gauss}
\liminf_{\ell \to \infty}
\frac 1\ell \log \bbP
\big( \tf_\ell (\gl, h; \go) \ge (1-\gd) \tf(\gl,0)
\big) \,\ge \, - \frac 12 h^2.
\end{equation}
To see this, we first observe that if $\bbP_{\ell,h}$ is the law of 
$(\go_1-h, \go_2 -h, \ldots, \go_\ell-h)$, then 
\begin{equation}
\label{eq:S}
\cS \left(\bbP_{\ell,h} \vert \bbP_{\ell,0}\right) \, :=\, 
\bbE_{\ell, h}\left[ \log \left( \frac{\dd \bbP_{\ell, h}}{\dd \bbP_{\ell,0}} \right)\right]\, =\, 
\frac 12 \ell h^2. 
\end{equation} 
We now recall the entropy inequality 
\begin{equation}
\label{eq:S-ineq}
\log  \left(\frac{\bbP_{\ell, 0} (E)}{\bbP_{\ell, h} (E)}\right)\, \ge \, 
- \frac 1{\bbP_{\ell, h} (E)} \left( \cS \left(\bbP_{\ell,h} \vert \bbP_{\ell,0}\right) + \frac 1 e
\right),
\end{equation}
which holds for arbitrary non-null events \cite[App.~A.2]{cf:Book},
and by choosing for $E$ the event in the left-hand side of
\eqref{eq:claim-Gauss} (call it $E_\ell$) we see that
\begin{equation}
\bbP_{\ell, h} \left(E_\ell\right) \, =\, \bbP_{\ell, 0}
\big( \tf_\ell (\gl, 0; \go) \ge (1-\gd) \tf(\gl,0)\big),
\end{equation}
where of course in the right-hand side we can write $\bbP$ instead of
$\bbP_{\ell,0}$.  But the existence of the infinite volume limit
\eqref{eq:fe}, together with the fact that $\tf(\gl,0)>0$ for $\gl>0$
(which follows from Proposition \ref{th:ublb}), guarantees that
$\lim_{\ell \to \infty }\bbP_{\ell, h} \left(E_\ell\right)=1$ and
this, combined of course with \eqref{eq:S} and \eqref{eq:S-ineq},
yields \eqref{eq:claim-Gauss}.

\begin{figure}
\begin{center}
\leavevmode
\epsfxsize =14 cm
\psfragscanon
\psfrag{0}[c]{$0$}
\psfrag{n}[c]{$n$}
\psfrag{Sn}[c]{$S_n$}
\psfrag{l}[c]{\small $\ell$}
\psfrag{b}[c]{\small $2\ell$}
\psfrag{c}[c]{\small $7\ell$}
\psfrag{d}[c]{\small $4\ell$}
\psfrag{e}[c]{\small $5\ell$}
\psfrag{f}[c]{\small $6\ell$}
\psfrag{g}[c]{\small $10\ell$}
\psfrag{h}[c]{\small $11\ell$}
\psfrag{i}[c]{\small $16\ell$}
\psfrag{m}[c]{\small $17\ell$}
\epsfbox{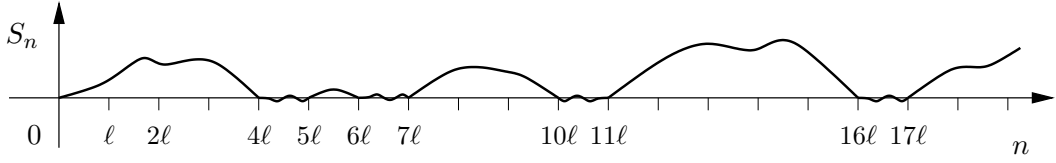}
\end{center}
\caption{\label{fig:3}
The novel strategy, compare with Figure~\ref{fig:2}, is based on targeting $q=0$ regions,
that is the (four) rich blocks are in this case the ones 
where $Q_\ell =o(\ell)$.
This is actually implemented in the proof  by a change
of measure argument. More precisely, the rich blocks in the Gaussian case are the
ones in which the charges look like the original charges $\go$
shifted down of $h$, so that $\go+h$ is a sequence of centered random
variables. Once again the lower bound is achieved by considering 
trajectories that stay in the upper half-plane outside of the rich blocks
and that touch the interface exactly at the beginning and
the end of a rich block. In a rich block, however, we keep the contribution
of all the trajectories, that can therefore oscillate between the two 
solvents in order to optimize the energetic gain. 
}
\end{figure}

The rest of the argument follows the line of \cite{cf:BG},
alternatively see \cite[Ch.~6]{cf:Book}, and it is based on chopping
the sequence of charges into portions of length $\ell$ and checking
whether the charge sub-sequence in the block is in $E_\ell$ (these are
the {\sl rich blocks}), namely whether $(\go_{j\ell +1}, \go_{j\ell +2}\ldots
 )\in E_\ell$ for $j=0, 1, \ldots$. These are of
course independent events giving rise to a Bernoulli sequence of
parameter $p(\ell):= \bbP\left( E_\ell\right)$.  Once $\go $ is
chosen, the rich blocks are identified and one estimates from below
$Z_{N, \go}$ by restricting to path configurations that visit the
sites $j\ell$ and $(j+1)\ell$ for all $j$'s for which the $j$-th block
is rich (and of course $(j+1)\ell\le N$) and that do not enter the
lower half plane outside of rich blocks (see Figure~\ref{fig:3} and
its caption for more details). The free energy bound one obtains is
\begin{equation}
\label{eq:Gcf}
\tf (\gl, h) \, \ge \, p(\ell) \left[
(1-\gd) \tf (\gl, 0) - (1+ \ga) \frac{h^2}2 + o_\ell(1) \right],
\end{equation}
therefore $\tf (\gl, h)>0$ if the term between square brackets is positive.  
Since $\ell$ can be chosen arbitrarily large and (then) $\gd$ arbitrarily small,
we obtain \eqref{eq:0strategy-Gauss}.
\qed 

\medskip

The result we have just stated becomes particularly effective
when coupled with the next statement.
\medskip

\begin{proposition}
\label{th:frombook}
Choose $K(\cdot)$ that satisfies \eqref{eq:K} and $\go_1 $
as in Definition~\ref{def:charges}. 
\begin{enumerate}
\item
If $\ga\ge 1$ then
\begin{equation}
\label{eq:frombook}
\lim_{\gl \searrow 0}
\frac 1{\gl^2} \tf (\gl, 0)\, =\, \frac 12.
\end{equation}
\item If $\ga \in (0,1)$ then 
for every  $\kappa>0$ we have
\begin{equation}
\label{eq:quasiexpl}
\liminf_{\gl \searrow 0}\frac 1 {\gl^2}\tf(\gl, 0)
\, \ge \
\frac{\kappa }{\Gamma(1-\ga)}
  \int_{0}^{\infty}
 \frac{\exp(-t)}{t^{1+\ga}} \bbE \left[ \log \cosh\left( z\sqrt{t/\kappa}\right)\right] \dd t
 \; -\; \kappa \frac{1-\ga}{\ga},
\end{equation}
where $z$ is a standard Gaussian random variable $\mathcal N(0,1)$ and 
$\Gamma(u) = \int_0^\infty dt  \; t^{u-1} {\rm e}^{-t}$.
\end{enumerate}
\end{proposition}
\medskip

\noindent
{\it Proof.} The case $\ga>1$ has been already considered in
\cite[Ch.~6]{cf:Book}, but we detail it here for completeness.
It is slightly more intuitive in this argument to work with
$\cZ_{N, \go}$ of \eqref{eq:cZ} (where it is understood that 
$\sign(S_n):=1-2\Delta_n$), and of course 
$\log(\cZ_{N, \go} / Z_{N, \go}) = \exp(\gl \sum_{n=1}^N \go_n)=\exp(o(N))$, since in this proof $h=0$.
The fact that, for every $\ga>0$,  $\limsup_{N\to \infty} (1/N) \bbE \log \cZ_{N, \go} 
\le \log \M (\gl)$ is an immediate consequence of Jensen's
inequality and of course $ \log \M (\gl)\sim \gl^2/2$ as $\gl \searrow 0$.
For what concerns the inferior limit, we still
apply Jensen's inequality, but not for switching
$\bbE $ and $\log$, rather $\log $ and $\bE$. Note that, with the 
notation $\go(j,k] := \sum_{n=j+1}^k \go_n$,
\begin{equation}
\label{eq:aflater}
\cZ_{N, \go} \,\ge \, 
\bE \left[
\prodtwo{j: \, j\ge 1}{\ \ \ \  \tau _j \le N}
\cosh
\left(
\gl \go (\tau_{j-1}, \tau_j]
\right)\right],
\end{equation}
which is obtained by integrating out the random signs $s=\{s_j\}_j $
(the inequality is due to neglecting the last, {\sl incomplete}, 
excursion, when present). With the notation $\psi(t):=\log \cosh(t)$, we get 
then
\begin{equation}
\label{eq:inv_annealed}
\frac 1N \bbE \log  \cZ_{N, \go} \,  \ge \, 
\frac 1N \bE\left[
\sumtwo{j: \, j\ge 1}{\ \ \ \ \tau _j \le N} \bbE \, \psi \left( \gl \go (\tau_{j-1}, \tau_{j}]\right) 
\right]
\stackrel{N\to \infty}\longrightarrow
\frac1{\bE[\tau_1]} \bE \bbE\psi \left( \gl \go (0, \tau_{1}]\right) ,
\end{equation}
where we have used also the law of large numbers.
Now we observe that, by the integrability properties of $\go_1$, we have that for every $n$ 
\begin{equation}
\label{eq:fatu1}
\lim_{\lambda \searrow 0}  \frac{1}{\gl ^2} \bbE\psi \left( \gl \go (0, n]\right)\, = \,
\frac{1}{2} \bbE \left[  (\go (0, n])^2 \right] \,= \,\frac{n}{2} \, .
\end{equation}
 Therefore, by applying Fatou's Lemma we infer from \eqref{eq:inv_annealed}
that
\begin{equation}
\label{eq:fatu2}
\liminf_{\gl \searrow 0} \frac{1}{\gl^2} \tf (\gl, 0) \,\ge \,  \frac1{\bE[\tau_1]} \sum_n \frac n2 K(n) \; =\; \frac 12,
\end{equation}
and the proof of the case $\ga>1$ is complete.

\smallskip

The case $\ga\le 1$ is more delicate and one has to go beyond the
direct use of Jensen's inequality. The estimates can be performed by
using a change-of-measure argument via a standard entropy inequality
(see {\sl e.g.} \cite[(A.10)]{cf:Book}). The case $\ga<1$ is fully
detailed in \cite[(6.17)]{cf:Book} (the proof in there follows from
sharpening an argument that appears in \cite{cf:BdH}). The case
$\ga=1$ is however not treated for the specific question we address
here, so we give a proof (which is close to the proof with which one
establishes \eqref{eq:quasiexpl}).

Set $\ga =1$. 
 Given $b>0$ define
$K_b(n)=K(n)\exp(-bn)/\left(\sum_{m} K(m)\exp(-bm)\right)$. 
Formula (6.13) in \cite{cf:Book} tells us that for every $b>0$
\begin{equation}
\label{eq:frombook3}
\tf(\gl,0) \, \ge \, \frac1{m_b} \sum_n K_b(n)
\bbE\left[ \psi (\gl \go(0,n]) \right] - \frac 1{m_b} \sum_n K_b(n) \log \left(
\frac{K_b(n)}{K(n)}
\right)\, =: \, e(b)-s(b),
\end{equation}  
where $m_b= \sum_n n K_b(n)$. 
In \cite[Proposition~B.2]{cf:Book} it is shown that $s(b)=o(b)$
(for  the case we are considering here, the result is a 
direct consequence of $\sum_n \exp(-bn )/n \sim
-\log b$, as $b \searrow 0$).
Moreover by using 
\begin{equation}
\label{eq:logcosh}
\log \cosh (x) \, \ge \, \frac 12 x^2 - \frac 1{12} x^4, \ \ \  \ x \in \R,
\end{equation}
and, by setting $b=\gl ^2$, from \eqref{eq:frombook3}
we obtain
\begin{multline}
\label{eq:ffrombook3}
\frac1{\gl^2}\tf(\gl,0) \, \ge \\
 \frac1{2 m_{\gl^2}} \sum_n 
K_{\gl^2}(n) \bbE\left[\left(\sum_{j=1}^n \go_j\right)^2\right]  - \frac {\gl^2}{12 m_{\gl^2}} 
\sum_n K_{\gl^2}(n) \bbE\left[\left(\sum_{j=1}^n \go_j\right)^4\right]  - \frac{s(\gl^2)}
 {\gl^2}
\\
\ge \,   \frac1{2 m_{\gl^2}} \sum_n 
K_{\gl^2}(n) n - \frac {c \gl^2}{m_{\gl^2}} 
\sum_n K_{\gl^2}(n) n^2  - \frac{s(\gl^2)}
 {\gl^2}\, \ge \, \frac 12 - 
 \frac {c^\prime \gl^2}{m_{\gl^2}} 
\sum_n \exp(-\gl^2 n)  - \frac{s(\gl^2)}
 {\gl^2}
 ,
\end{multline}
where $c$ is a constant that depends only on the fourth moment of 
$\go_1$ and $c^\prime$ comes from approximating 
$K(n)$ with its limit behavior and from neglecting $\sum_n K(n) \exp(-\gl^2 n)$
in the denominator (say, for $\gl \le 1$). Therefore 
$\liminf_{\gl \searrow 0 } \tf (\gl, 0)/ \gl^2 \ge 1/2$ since
$m_{\gl^2} $ diverges as $\gl \searrow 0$ and $s(\gl^2)/ \gl^2$ tends to zero,
as pointed out before. 
\qed

\medskip

\begin{rem}\label{rem:nienteunmezzo}
  \rm For $\ga=1/2$, from a numerical estimation of the right-hand
  side of \eqref{eq:quasiexpl}, we obtain that $l:=\liminf_{\gl
    \searrow 0} \tf(\gl,0) /\gl^2 > 0.227$. On the other hand, one
  would need $l>1/3$ for our new strategy to be better than the older
  one, i.e., to be able to prove that
  $\liminf_{\gl\searrow0}h_c(\gl)/\gl>2/3$ with our method (just
  recall \eqref{eq:0strategy-Gauss}). An evaluation of $l$ by using
  the transfer matrix method (with the software developed in
  \cite{cf:CGG}) for small values of $\gl$ suggests that $l$ is below
  $1/3$, even if it looks rather close to it, which in particular
  tells us that the neutral-stretch strategy is probably better than
  the old one down to $\ga$ very close to $1/2$.
\end{rem}

\medskip

\noindent
{\it Proof of Theorem~\ref{th:neutral}, Gaussian charges.}
From Propositions \ref{th:0strategy-Gauss} and \ref{th:frombook}(1)
one gets immediately the bound $\liminf_{\gl\searrow0}h_c(\gl)/\gl\ge 1/\sqrt{1+\ga}$ of \eqref{eq:neutral1}.
The fact that the same (inferior) limit is  not smaller than $1/2$ for $\ga\ge1$ follows since for $\ga>1$
\begin{equation}
\begin{split}
  \frac1N\bbE\log   Z_{N,\go}\, &\ge\, -h\gl+\frac1N\bE
  \sumtwo{j: \, j\ge 1}{\ \ \ \ \tau _j \le N} \bbE \left[
    \psi \left( \gl h(\tau_j-\tau_{j-1})+\gl \go (\tau_{j-1}, \tau_{j}]\right) 
  \right]\\
  &\stackrel{N\to \infty}\longrightarrow\, 
-h\gl+\frac1{\bE[\tau_1]} \bE \bbE\psi \left(h\gl\tau_1+ \gl \go (0, \tau_{1}]\right) ,
\end{split}
\end{equation}
(which is analogous to \eqref{eq:inv_annealed}) and from the fact that
\eqref{eq:fatu1} still holds with $\psi(\gl\go(0,n])$ replaced by
$\psi(\gl h n+\gl\go(0,n])$, if $ h=O(\gl)$.

Formula \eqref{eq:neutral2}, 
for $\ga<1$ but close to $1$,
follows from  Propositions \ref{th:0strategy-Gauss}
and \ref{th:frombook}(2), plus the observation that
the limit of the right-hand side of \eqref{eq:quasiexpl} as $\ga\nearrow 1$
is equal to $1/2$ (this has been already pointed out and detailed in 
\cite[Remark~6.4]{cf:Book}). 
A more quantitative bound can be obtained as follows.
By using the inequality \eqref{eq:logcosh},
one can bound from below the right-hand side in \eqref{eq:quasiexpl}
by a quantity that can be explicitly computed: 
\begin{equation}
\label{eq:quasiexpl2}
\liminf_{\gl \searrow 0}\frac 1 {\gl^2}\tf(\gl, 0)
\, \ge \frac 12 -\frac{1 -\ga}{4\kappa}- \kappa \left(\frac{1-\ga}\ga\right)
\, \stackrel{\kappa=\sqrt{\ga}/2}= 
\frac 12- \frac{1-\ga}{\sqrt{\ga}} \, ,
\end{equation}
where we used that $\Gamma(2-\alpha) = ( 1- \alpha) \Gamma(1-\alpha)$.
From \eqref{eq:0strategy-Gauss} we see that
\eqref{eq:neutral2} holds
if $\liminf_{\gl \searrow 0}  \tf(\gl,0)/\gl^2 > 1/(2(1+\ga))$, so that
from \eqref{eq:quasiexpl2} we have that
\eqref{eq:neutral2} holds
if $\ga >  0.800981\ldots$.
A numerical evaluation of the full expression
\eqref{eq:quasiexpl} shows that \eqref{eq:neutral2}
holds down to $\ga =0.65$.
\qed

\medskip

\noindent
{\it Proof of Theorem~\ref{th:neutral}, general charges.}
As explained in section \ref{sec:rev}, one
expects a universal behavior for  $\lambda\searrow 0$ and the Gaussian
bounds should remain in force. And in fact
the strategy that leads to Proposition~\ref{th:0strategy-Gauss} can be
generalized without much effort, but Proposition~\ref{th:0strategy-Gauss} becomes 
more complex to state and somewhat involved.  It is therefore preferable 
to approach the problem under a slightly different angle, which in the end 
is just dealing
with Proposition~\ref{th:0strategy-Gauss} and Proposition~\ref{th:frombook}
at the same time. 

The key point is to replace the change of measure used in the Gaussian
case with the standard {\sl tilting} procedure of Cram\`er Large
Deviation Theorem (which coincides with a shift of the mean for
Gaussian variables).  We therefore modify the law of $(\go_1, \ldots,
\go_\ell)$ by introducing the relative density $f_\ell
(\go):=\prod_{j=1}^\ell \exp( \mu_h \go_j)/ \M ( \mu_h)^\ell$, and $\mu_h$
is chosen so that $\bbE[ f_\ell (\go) \go_1]=- h$.  The target event
$E_\ell$ in \eqref{eq:claim-Gauss} should now be replaced by
\begin{equation}
E_\ell \, =\, \left\{ \go:\, \tf _\ell (\gl, h; \go) \ge (1-\gd) \tilde \tf_h(\gl)\right\},
\end{equation}
where $\tilde  \tf_h(\gl)$ is the free energy of the model 
with IID charges $\tilde \go$, with the law of $\tilde \go_1$
given by the law of $\go_1$ times the density $  \exp( \mu_h \go_1)/ \M ( \mu_h)$.
Note that the new charges $\tilde \go$ have mean $-h$ and, in general, their variance is not equal to one, but it converges to 1 as $h$ vanishes.
This induces a number of changes that lead to
replacing \eqref{eq:Gcf} by
\begin{equation}
\tf (\gl, h) \, \ge \, p(\ell) 
\left[ (1-\gd) \tilde  \tf _h (\gl) - (1+\ga) \gS (-h) + o_\ell (1)\right],
\end{equation}
where $\gS(\cdot)$ is the Cramer functional of the law of $\go_1$,
{\sl i.e.} the Legendre transform of $\M(\cdot)$. It is well known
that $\gS (-h) \sim h^2/2$ as $h \searrow 0$, thus for every
$\epsilon>0$ there exists $h_\epsilon>0$ such that $(\gl, h) \in \cL$
if
$h< h_\epsilon$ 
\and
 \begin{equation}
 (1-\gd) \tilde\tf _h (\gl) \, \ge\,  (1+\ga)(1+\epsilon) \frac{h^2}{2}.
 \end{equation}

 At this point, $\tilde \tf _h (\gl)$ can be bounded from below
 precisely as it is done in Proposition~\ref{th:frombook}. Choosing $h
 = m \lambda$, the argument in that case always leads to estimating
 the moments (in fact, the second and the fourth moments suffice) of
 $\go_1$, that has to be replaced with the centered variable $\tilde
 \go_1 -h$ of variance $1+O(\lambda)$ (and $\bbE[(\tilde \go _1+h)
 ^4]=\bbE[ \go _1 ^4] (1+O(\lambda)$).

This completes the proof of Theorem~\ref{th:neutral}.
\qed

\section{Fractional moments and upper bounds}
\label{sec:frac}

In this section we prove Theorem \ref{th:ub}.

We use the short-cut notation
$Z_N:= Z_{N, \go}^\rc$ and $Z_{a,b}:=Z^\rc_{(b-a),\theta^a\go}$
for $a<b$, where $\theta$ is the shift operator such that 
$(\theta\go)_n=\go_{n+1}$.
We start by pointing out that, by integrating out
the $\{ s_j \}_j$ variables,  we can write
\begin{equation}
Z_N\, =\, 
\bE\left[\prod_{j: \tau_j \le N}
\gp\big(\gl\go(\tau_{j-1},\tau_j]+\gl h(\tau_j-\tau_{j-1})\big); \, N\in\tau\right],
\end{equation}
with
$\gp(t)=(1+\exp(-2t))/2 $ and $\go(j,k] := \sum_{n=j+1}^k \go_n$.

One of the two key ingredients for the proof is the following
decomposition of the partition function, which is just based on
partitioning the space of trajectories according to the location of
the first point of $\tau$ after $k$ (call it $j$) and the last point
before $k$  (call it $i$):
\begin{equation}
\label{eq:decomp}
Z_N=\sum_{j=k}^N Z_{j,N}
\sum_{i=0}^{k-1}K(j-i)\gp\big(\gl\go(N-j,N-i]+\gl h(j-i)\big)
Z_{i}.
\end{equation}
The second key ingredient is the use of fractional moments, so we set
for $0<\gamma<1$
\begin{equation}
A_N\, :=\, \bbE\left[(Z_N)^{\gga}\right],
\end{equation}
(we do not make the $\gamma$-dependence explicit in the notation).
From \eqref{eq:decomp} and the basic inequality 
\begin{equation}
\label{eq:basicineq}
\left(\sum a_i\right)^\gga \le \sum a_i^\gga, 
\end{equation}
that holds whenever $a_i\ge 0$ for every $i$ if $0<\gamma<1$,
we obtain
\begin{equation}
\label{eq:imp}
A_N\, \le\,  \sum_{j=k}^N A_{N-j}\sum_{i=0}^{k-1}B(j-i)A_i,
\end{equation}
where
\begin{equation}
B(j)\, :=\, K(j)^\gga\bbE\left[\Big(\gp\big(\gl\go(0,j]+\gl hj\big)\Big)^\gga\right].
\end{equation}
For later use we point out that, again thanks to \eqref{eq:basicineq},
\begin{equation}
\label{eq:bbd}
B(j)\, \le\,  K(j)^{\gga}2^{-\gamma}
\left[\exp\Big( j\big(\log M(-2\gga\gl)-2\gga\gl h \big)\Big)+1\right].
\end{equation}

\medskip

\begin{lemma}
\label{th:lem1}
If there exist $\gamma\in (0,1)$ and $k\in \N$ such that
\begin{equation} 
\label{eq:adem}
\sum_{j=k}^{\infty}\sum_{i=0}^{k-1}B(j-i)A_i\le 1,
\end{equation}
then $\sup_N A_N < \infty$.
\end{lemma}

\medskip

\noindent
{\it Proof.} Combining  \eqref{eq:imp} and the hypothesis \eqref{eq:adem}
one readily obtains for every $N\ge k$
\begin{equation}
\label{eq:newiterineq}
A_N \, \le \, \max_{j=0, \ldots, N-k} A_j.
\end{equation}
If one sets $A^\star_k:= \max_{j=0, \ldots, k-1} A_j$, then from
\eqref{eq:newiterineq} one has $ A_N \le A^\star_k$ for every $N$. 
\qed

\medskip

\begin{rem}
\label{rem:deloc}
\rm
Note that if $\sup_N A_N < \infty$ one has
\begin{equation}
\lim_{N \to \infty} \frac{1}{N}\bbE\left[\log Z_N\right]\, \le\, \lim_{N \to \infty} \frac{1}{\gga N} \log A_N=0
\end{equation}
and therefore $\tf (\gl,h)=0$.
\end{rem}

\subsection{Gaussian disorder and $\ga>1$}

\label{sec:a>1g}

Let us set $\alpha>1$ and $\go_1\sim \cN (0,1)$, and fix
$\gl_0<\infty$. We are going to show that if $\rho<1$ is sufficiently
close to one, $\tf(\gl,\rho\gl)=0$ for $\gl\le \gl_0$.  To
this purpose we fix $\gga<1$ such that $\gga(1+\alpha)>2$ and notice
that
\begin{equation}
\label{eq:Gauss-1}
\log \M(-2\gamma\gl)-2\gamma\rho\gl^2=2\gamma\gl^2(\gamma-\rho)<0,
\end{equation}
provided that $\rho>\gamma$. Therefore, thanks also 
to \eqref{eq:adem} and
\eqref{eq:bbd}, 
it is sufficient to show that
\begin{equation}
\label{eq:tbsub}
\sum_{j=k}^\infty\sum_{i=0}^{k-1}   K(j-i)^{\gga}A_i\le 2^{\gga -1}.
\end{equation}
Computing the sum over $j$, we see that there exists $C>0$ depending on
$K(\cdot)$ such that the 
left-hand side in \eqref{eq:tbsub} is smaller than
\begin{equation}
\label{eq:truc}
C
\sum_{i=0}^{k-1} \frac{A_i}{(k-i)^{\gga(\alpha+1)-1}}. 
\end{equation}
In view of the last two formulas we see that it is sufficient to show
that there exist $\rho\in (\gamma,1)$ 
such that for 
$\gl\le \gl_0$ and  $h=\rho\gl$ it is possible to find $k\in \N $ such that 
\begin{equation}
\label{eq:tobeshown}
\sum_{i=0}^{k-1} \frac{A_i }{(k-i)^{\gga(\alpha+1)-1}}\, \le\, c_1 := 
\frac{2^{\gamma-1}}C. 
\end{equation}
We choose $k$ to be equal to the integer part of  $ 1/(\gl^2(1-\rho))$
(but there is no loss of generality in choosing $ 1/(\gl^2(1-\rho))\in \N$, so we do that)
and we note that $k$ is large (uniformly in 
$\gl\le\gl_0$) if $\rho$ is close to $1$. Hence, by Jensen's inequality we have
\begin{equation}
\label{eq:Gauss+1}
A_i \le [\bbE\,Z_i]^\gamma \, \le\,  \exp\big(2i(1-\rho)\gl^2 \gamma \big)\, \le\,  e^{2 \gamma},
\end{equation}
for $i\le k$. Therefore, if $R$ is an  integer number smaller than $k$
\begin{equation}
\label{eq:dkeg}
\sum_{i=0}^{k-R}\frac{A_i }{(k-i)^{\gga(\alpha+1)-1}}\, \le\,  e^{2\gga}\sum_{l=R}^\infty \frac{1}{l^{\gga(\alpha+1)-1}}.
\end{equation}
The right-hand side can be made arbitrarily small by choosing $R$ large
(of course this requires $k$ to be large, but this is again the requirement of choosing
 $\rho$  close to $1$), because
$\gga$ is such that $\gga(\alpha+1)>2$.  We choose $R$ large enough so
that such an expression is smaller that $c_1/2$.  It suffices now to
prove that
\begin{equation}
\sum_{i=k-R+1}^{k-1}\frac{A_i}{(k-i)^{\gga(\alpha+1)-1}}\, \le \, c_1/2,
\end{equation}
for  $\rho$ sufficiently close to $1$.
In analogy with \cite{cf:DGLT}, we are going to shift the
random variables $\go_i$ in order to bound $A_i$. We define
the shifted charges by introducing the new law $\bbP _{N,y}$
\begin{equation}
\label{eq:shift-1}
  \frac{\dd \bbP_{N,y}}{\dd \bbP}(\go)=\exp\left(y\sum_{i=1}^N \go_i -Ny^2/2\right),
\end{equation}
and control $A_i$ by applying the H\"older inequality with $p=1/\gamma$
and $q=1/(1-\gamma)$:
\begin{equation}
\label{eq:Hold}
A_N\, =\, 
\bbE_{N,y}\left[\left(Z_{N}\right)^\gamma\frac{\dd \bbP}
  {\dd \bbP_{N,y}}(\go)\right]
\, \le \, 
\left(\bbE_{N,y}\left[Z_{N}\right]\right)^\gga
\exp\left(\frac{\gga}{2(1-\gga)}y^2N\right). 
\end{equation}
We apply \eqref{eq:Hold} with $y=\gl(1-\rho)^{1/2}=1/\sqrt k$, $N=i$
and $k-R<i < k$.  Note that with this choice the exponential term in
the right-most side of \eqref{eq:Hold} is bounded by the constant
$\exp\big(\gga/(2(1-\gga))\big)$, that $\bbE_{i,y}\left[Z_i\right]$
coincides with
\begin{equation}
\label{eq:expr-ub}
\bE\left[
\exp\left(
2\gl ^2\left( (1-\rho)- (1-\rho)^{1/2}\right) \sum_{n=1}^{i} \gD_n
\right); \,  i\in\tau
\right]
\end{equation}
and that for $\rho$ close to $1$
\begin{equation}
\label{eq:Gauss+2}
2\gl^2 \left( (1-\rho)- (1-\rho)^{1/2}\right)\, \le \, -\gl ^2(1-\rho)^{1/2}.
\end{equation}
Remarking that 
\begin{equation}
\min_{k-R<i < k}  i\gl^2 (1-\rho) \ge 1/2,
\end{equation}
 for $k$ large, we  see
that the quantity in \eqref{eq:expr-ub} 
can be made arbitrarily small, uniformly in $\gl\le \gl_0$ and 
$k-R<i<k$ by choosing $\rho$ suitably  close to $1$,
because
\begin{equation}
\label{eq:expr-ub2}
 \lim_{N \to \infty}
\bE\left[
\exp\left(- \frac q N
 \sum_{n=1}^{N} \gD_n
\right)
\right]\, =\, \exp(-q/2),
\end{equation}
which follows from the Dominated Convergence Theorem
since $(1/N) \sum_{n=1}^{N} \gD_n$ tends  almost surely to $1/2$.
The latter statement follows by observing that if we set
$Y_N:= \max\{n:\tau_n\le N\}$ (number of renewals up to $N$),
by the law of large numbers we have that $Y_N/N$ tends a.s.
to $1/ \bE[\tau_1]$ as $N $ tends to infinity and
\begin{equation}
\label{eq:uselln}
\frac 1N \sum_{n=1}^N \gD_n \, \ge \, 
\frac {Y_N}N \cdot\frac 1{Y_N}\sum_{j=1}^{Y_N} (\tau_j-\tau_{j-1}) \ind_{s_j=-1}.
\end{equation}
Since 
 $\{ (\tau_j-\tau_{j-1}) \ind_{s_j=-1}\}_{j=1, 2, \ldots}$ is an IID sequence and 
$\bE[ \tau_1 \ind_{s_1=-1}]=\bE[\tau_1]/2$, again by the law of large numbers
 the right-hand side in  \eqref{eq:uselln} converges almost surely
to $1/2$. We can reverse the inequality in \eqref{eq:uselln} by 
summing over $j$ up to $Y_N+1$, so that \eqref{eq:expr-ub2}
is proven.

On the other hand 
\begin{equation}
  \sum_{i=k-R+1}^{k-1}\frac{A_i}{(k-i)^{\gga(\alpha+1)-1}}\,\le\, 
  \max_{k-R < j < k} A_j \, 
\sum_{i=1}^{\infty}\frac{1}{i^{\gga(\alpha+1)-1}},
\end{equation}
so that, since the sum converges, the right-hand side can be made
arbitrarily small, hence smaller than $c_1/2$, and we are done.

The fact that $h_c(\gl)<h^{(1)}(\gl)$ for every $\gl>0$ is a direct consequence
of the fact that  $\gl_0$ is arbitrary.
\qed

\medskip

\subsection{Gaussian disorder and  $\alpha<1$}

Again, fix $\gl_0<\infty$.  We let $h=h(\gl)= \gl\left(1 -
  \frac{c}{|\log c\gl^2|}\right)$ and we aim at proving \eqref{eq:adem}
for some $c>0$, uniformly in $0<\gl\le \gl_0$.  We choose $k$ to be
equal to (the integer part of) $\frac{|\log c\gl^2|}{c\gl^2}$ and
$\gga=1-(\log k)^{-1}$. Notice that $k$ can be made arbitrarily large
(and therefore $\gamma$ close to $1$) uniformly for all $\gl\le \gl_0$
choosing $c$  small enough. Therefore, we will consider that $k$ is large
when we need to.

We use \eqref{eq:bbd} to find a simple bound on $B(\cdot)$. Since
\begin{equation}
\label{eq:Gauss-2}
  \log(M(-2\gga\gl))-2 \gga \gl h(\gl)\, =\, 2\gga\gl^2\left(\gga -
    \frac{h(\gl)}{\gl}\right)\le 0,
\end{equation}
when $c$ is well chosen,
we have
\begin{equation}
B(j)\le K(j)^{\gga} 2^{1-\gga}.
\end{equation}
Therefore the condition \eqref{eq:adem} will be fulfilled if we can show that
\begin{align}
  \sum_{j=k}^{\infty}\sum_{i=0}^{k-1}\frac{A_i}{(j-i)^{(\alpha+1)\gga}}\le
  c_2 \label{eq:inf2}
\end{align}
for a suitable constant $c_2$ which is independent of $c$ or $\gl$.
First of all we get rid of $\gga$ in the denominator as follows:
\begin{align}
  &\sum_{j=k}^{\infty}\frac{1}{(j-i)^{(\alpha+1)\gga}}=
\sum_{j=(k-i)}^\infty j^{-(\alpha+1)\gga}\\
  &\quad \le \sum_{j=(k-i)}^{k^6}
  j^{-(\alpha+1)}\exp\left(\frac{(\alpha+1)\log j}{\log
      k}\right)+\sum_{j= k^6+1}^\infty j^{-(1+\frac{\alpha}{2})}\le c_3
  (k-i)^{-\alpha},
\end{align}
for some constant $c_3<\infty$, provided $c$ is such that
$(\alpha+1)\gga\ge 1+(\alpha/2)$. Hence \eqref{eq:inf2} will be
satisfied if 
\begin{align}
\sum_{i=0}^{k-1}\frac{A_i}{(k-i)^{\alpha}}\le \frac{c_2}{c_3}. \label{eq:inf1}
\end{align}
To estimate $A_i$ we use \eqref{eq:Hold} with
$y={\sqrt{c}\gl}/{|\log c\gl^2|}$, $N=i$ and $i<k$.
With these settings, one can check that the exponential term 
in the right-hand side of \eqref{eq:Hold} is bounded
by a constant $c_4$ that does not depend on $\gl$ or $c$ (and therefore will 
be harmless), and that $\bbE_{i,y}\left[Z_{i}\right]$ coincides
with
\begin{align}
  \bE\left[\exp\left(2\frac{\gl^2}{|\log
        c\gl^2|}\left(c -\sqrt{c}\right)\sum_{n=1}^i
      \gD_n\right);i\in\tau\right]\le 
\bE\left[\exp\left(-\frac i{\sqrt c k}\frac{\sum_{n=1}^i
      \gD_n}i\right);i\in\tau\right]
 \label{eq:ai}
\end{align}
for $c$ small enough.
The last expression is in any case smaller than $\bP[i\in\tau]$, 
which is itself bounded above by $c_5
i^{\alpha-1}$ for some constant $c_5$ depending on $K(\cdot)$ (see
\cite[Theorem B] {cf:Doney}).

We need a better upper bound for large $i$. This is provided by
\smallskip

\begin{lemma}
\label{th:Deltas}
Assume that $0<\ga\le 1$. Then,
  \begin{eqnarray}
\label{eq:Deltas}
    \lim_{q\to\infty}\limsup_{N\to\infty}
\bE\left[\left.\exp\left(-\frac qN{\sum_{n=1}^N
      \gD_n}\right)\right|N\in\tau\right]=0.
  \end{eqnarray}
\end{lemma}
\medskip

We fix a small $a>0$, and consider $ ak+1\le i<k$.
Since therefore $i/k>a$, thanks to Lemma \ref{th:Deltas} one deduces that,
if one chooses $c$ sufficiently small (how small, depending on $a$),
the quantity  in the right-hand side of \eqref{eq:ai} is bounded above by 
$a \bP[i\in\tau]$.

We can summarize our result concerning $A_i$ as follows:
\begin{equation}
A_i\, \le \, 
\begin{cases}
\label{eq:uffa}
   \left(\bP(i\in \tau)\right)^{\gga}\le c_6 i^{\alpha-1}  &
\text{ for } i\le ak,\\
   \left(a \bP(i\in \tau)\right)^{\gga}\le c_6 a^\gga
  i^{\alpha-1}\le c_6 a^{\gga+\alpha-1} k^{\alpha-1} 
   &\text{ for
  } ak+1 \le i \le k-1,
\end{cases}
\end{equation}
where in both cases we used that $i^{1-\gamma}\le c_7$ for $i\le k$.
Hence
\begin{equation}
\label{eq:hence}
\sum_{i=0}^{ak}\frac{A_i}{(k-i)^\alpha}\le c_8 a^{\alpha}\  \  \text{ and } \ \
\sum_{i=ak+1}^{k-1}\frac{A_i}{(k-i)^\alpha}\le c_9 a^{\gga+\alpha-1},
\end{equation}
where we just used the previous inequalities to estimate $A_i$. Since $a$
can be chosen arbitrarily small, \eqref{eq:inf1} is satisfied.

\qed

{\sl Proof of Lemma \ref{th:Deltas} for $0<\ga<1$}.
First of all, we claim that it is sufficient to prove \eqref{eq:Deltas}
for the unconditioned measure, i.e., that
\begin{eqnarray}
\label{eq:Deltas1}
    \lim_{q\to+\infty}\limsup_{N\to\infty}
\bE\left[\exp\left(-q\frac{\sum_{n=1}^N
      \gD_n}N\right)\right]=0.
  \end{eqnarray}
 Indeed, with $X_N :=\max\{n=0,1, \ldots, N/2:\, n \in \tau\}$
(last renewal epoch up to $N/2$), one has 
\begin{eqnarray}
&&  \bE\left[\left.\exp\left(-q\frac{\sum_{n=1}^N
      \gD_n}N\right)\right|N\in\tau\right]\le
\bE\left[\left.\exp\left(-q\frac{\sum_{n=1}^{N/2}
      \gD_n}N\right)\right|N\in\tau\right]\\\nonumber
&&= \sum_{k=0}^{N/2} \bE 
\left[\left.\exp\left(-q\frac{\sum_{n=1}^{N/2}
      \gD_n}N\right)\right|  \,X_N=k
\right] \bP\left( X_N=k \big \vert \, N \in \tau\right)
\end{eqnarray}
where we used the renewal property in the second step.
Next, it is not difficult to see that for every $k =0, 1, \ldots, N/2$
\begin{equation}
\bP\left( X_N=k \big \vert \, N \in \tau\right) \le c_{10} 
\bP\left( X_N=k\right)
\end{equation}
(this is detailed for instance in the proof of \cite[Lemma 4.1]{cf:DGLT}) and 
the claim follows.

To show \eqref{eq:Deltas1}, 
note that, again with the notation $Y_N:=\max\{n:\tau_n\le N\}$, for every $\gep\in(0,1)$
  we have
\begin{equation}
\bP\left( \frac{Y_N}{N^\ga} \le \gep \right) \, =\, 
\bP\left( \tau_{\gep N^\ga}\ge N\right)\, =\, 
\bP\left( \frac{\tau_{\gep N^\ga}}{(\gep N^\ga)^{1/\ga}}\ge \gep^{-1/\ga}\right)\stackrel{N \to \infty}\longrightarrow G_\ga (\gep^{-1/\ga}),
\end{equation}
where we have assumed $\gep N^\ga\in \N$  and $G_\ga(\cdot)$ is the integrated tail probability function
of an $\ga$-stable variable, see 
 \cite[Th.s 1 and 2,
pp.~448-449]{cf:Feller2} 
from which 
we extract also that $G_\ga (\gep^{-1/\ga}) \stackrel{\gep \searrow 0}
\sim \gep C_K /\Gamma (1-\ga)$ (recall \eqref{eq:K}).  Therefore 
\begin{equation}
\bP\left( \frac{Y_N}{N^\ga} \le \gep \right) \,\le \,2 C_K 
\gep/\Gamma(1-\ga),
\end{equation}
for $N$ sufficiently large and $\gep \in (0, 1)$. On the other hand 
\begin{equation}
\begin{split}
\bE \left[
\exp \left(
-\frac q N \sum_{n=1}^{N}\gD_n\right)\right]\, &\le \,
\bE\left[
\prod_{i=1}^{Y_N} \left(
\frac{\exp(-q (\tau_{i}-\tau_{i-1})/N)+1}{2}
\right)
\right] 
\\
&\le \, \bE \left[
\prod_{i=1}^{\gep N^\ga} \left(
\frac{\exp(-q (\tau_{i}-\tau_{i-1})/N)+1}{2}
\right)
\right] \, +
\, \bP\left( \frac{Y_N}{N^\ga} \le \gep \right)
\\
&\le \, \bE\left[\frac{\exp(-q\tau_1/N)+1}2\right]^{\gep N ^\ga} + 
2 C_K 
\gep/\Gamma(1-\ga).
\end{split}
\end{equation}
Observe that 
\begin{equation}
1-\bE\left[\frac{\exp(-q\tau_1/N)+1}2\right] 
\stackrel{N \to \infty}\sim
\frac{C_K q^\ga}{2\ga N^\ga} \Gamma (1-\ga), 
\end{equation}
so we have that for some constant $c_{11}(\ga)$ (which depends only on $\ga$)
\begin{equation}
\limsup _{N \to \infty } \bE \left[
\exp \left(
-\frac q N \sum_{n=1}^N\gD_n\right)\right] \, \le \, 
\exp\left(-c_{11}(\ga) (C_K  \gep) q^\ga \right) \, +\,  \frac{2 C_K 
\gep}{\Gamma(1-\ga)} \, \le \, 4\frac{(\log q)^2}{q^\ga\,\Gamma(1-\ga)},
\end{equation}
where in the last step we have chosen
$  C_K\gep= \log (q)^{2}/ q^\ga$ 
and we assumed that $q\ge q_0(\ga)$ with $q_0(\ga)$ is sufficiently large. 
The proof of Lemma \ref{th:Deltas} for $0<\ga<1$
is therefore complete.
 \qed

\subsection{Gaussian disorder and  $\alpha=1$}
The proof is very similar to that of the case $\ga<1$, and therefore
we point out only the necessary modifications.  No changes are
needed up to formula \eqref{eq:ai} and, again, we let $a$ be a small
positive number. To avoid repetitions, it will be understood that $c$
is chosen sufficiently small (how small depends on $a$), so that
$1/k$ and $1-\gamma=1/\log k$ can be made arbitrarily small with $a$
fixed.  Using simply the fact that $A_i\le c_{12}$, we obtain in analogy
with the first bound of \eqref{eq:hence}
\begin{eqnarray}
  \sum_{i=0}^{ak}\frac{A_i}{(k-i)}\le 2c_{12} a.
\end{eqnarray}
As for the values $i> ak$, we use the fact that (see \cite[Th. A.6]{cf:Book})
\begin{eqnarray}
\label{eq:A6}
  \bP(N\in\tau)\stackrel{N\to\infty}\sim \frac{c_{13}}{\log N}
\end{eqnarray}
and we claim that Lemma \ref{th:Deltas} still holds for $\alpha=1$
(this will be proven in a while) so that, in analogy with \eqref{eq:uffa},
$A_i\le a^\gamma\bP(i\in\tau)^\gamma$ for all $ak<i<k$.
Then, via \eqref{eq:A6} and recalling that $\gamma=1-1/\log k$,
\begin{eqnarray}
\sum_{i=ak+1}^{k-1}\frac{A_i}{(k-i)}\le c_{14}\log k \max_{ak<i<k}A_i\le
c_{15} a^\gamma
\end{eqnarray}
where we used also the fact that, choosing $c$ small enough, we can assume
$\log k/(\log ak)<2$.
This concludes the proof since \eqref{eq:inf1} is satisfied if $a$ is small.

\qed

\medskip

{\sl Proof of Lemma \ref{th:Deltas} for $\ga=1$}.  Again, it is
sufficient to prove the claim for the unconditioned measure, i.e., to
show \eqref{eq:Deltas1}.  
If we define
the event
\begin{eqnarray}
  \label{eq:AN}
  E_N:=\left\{(\tau_i-\tau_{i-1})\le \frac N{\sqrt{\log N}}\;\mbox{for every}\;
i\le Y_N+1\right\},
\end{eqnarray}
we have
\begin{eqnarray}
\label{eq:Deltas2}
\bE\left[\exp\left(-q\frac{\sum_{n=1}^N
      \gD_n}N\right)\right]&\le&
 \bE\left[\prod_{i\le Y_N}\left(\frac{1+e^{-q(\tau_i-\tau_{i-1})/N}}
2\right)\right]\\\nonumber
&\le&
      \bE\left[e^{-\frac q{4N}\tau_{Y_N}}\ind_{\{E_N\}}\right]+\bP\left(E_N^c\right),
  \end{eqnarray}
for $N$ sufficiently large
(we simply used that 
$\max_{i\le Y_N}(\tau_i-\tau_{i-1})/N$ 
tends to $0$ for $N\to\infty$ if $E_N$ is realized). Note also that, by the definition of $E_N$, 
$\tau_{Y_N}/N\ge 1/2$ for $N$ large if $E_N$ is realized.
Therefore, 
\begin{eqnarray}
 \limsup_{N\to\infty}\bE\left[\exp\left(-q\frac{\sum_{n=1}^N
      \gD_n}N\right)\right] \le e^{-\frac q{8}}+\limsup_{N\to\infty}
\bP\left(E_N^c\right).
\end{eqnarray}
To show that the probability of $E_N^c$ tends to zero with $N$, 
we start by observing that 
from the tail behavior of $K(\cdot)$ it follows that
\begin{eqnarray}
\label{eq:AN1}
\bP\left( \text{there exists } i\le \frac{N}{(\log N)^{3/4}}
\text{ such that} (\tau_i-\tau_{i-1})\ge 
  \frac N{\sqrt {\log N}}\right)\le \frac{c_{16}}{(\log N)^{1/4}}.
\end{eqnarray}
On the other hand,
\begin{equation}
\label{eq:AN2}
\lim_{N \to\infty}  \bP\left(Y_N\ge \frac{N}{(\log N)^{3/4}}\right)\le
\frac{(\log N)^{3/4}}N  \bE[Y_N]=\frac{(\log N)^{3/4}}N
\sum_{i=1}^N \bP(i\in\tau)
\le \frac{c_{17}}{(\log N)^{1/4}},
\end{equation}
where we used \eqref{eq:A6}. 
From \eqref{eq:AN1}-\eqref{eq:AN2} we directly see that
$\lim_{N\to\infty} \bP(E_N^c)=0$ and the proof is complete.

\qed

\subsection{Proof of Theorem~\ref{th:ub}: the general case}
Once again, going beyond the case of Gaussian disorder
requires only a bit of care in some steps. In both case ($\alpha>1$ and $\alpha\le 1$), we have to prove
\begin{align}
h_c(\gl)\le h^{(1)}(\rho\gl)\quad \text{ for every }  \ \gl\le \gl_0,
\end{align}
where $\rho$ is a fixed constant, chosen close to $1$, in the case
$\alpha>1$ (we can get the slope result from it then), and
$\rho=\rho(\gl)=1-\frac{c}{|\log c \gl^2|}$ for a small $c$ for the
case $\alpha\le 1$. We take here a quick look at what needs to be
changed from the proof of the Gaussian case.

Precisely \eqref{eq:Gauss-1}, \eqref{eq:Gauss-2}
have to be replaced by the observation that
\begin{equation}
\log \M (-2 \gamma \gl) - 2 \gamma \gl h^{(1)}(\rho \gl)\, =\, 
2 \gamma \gl \left[ h^{(1)}(\gamma \gl)-h^{(1)}(\rho \gl)\right],
\end{equation}
and that the right-hand side is negative if $\rho>\gamma$ 
because $h^{(1)}(\cdot)$ is increasing. 
For what concerns instead
\eqref{eq:Gauss+1} (and the same bound has to be used for \eqref{eq:expr-ub} and \eqref{eq:ai}) the bound one has to use are 
\begin{equation}
\label{eq:nonGauss-2}
\log \M (-2  \gl) - 2 \gl h^{(1)}(\rho \gl)\, =\, 
\log \M (-2  \gl)- \frac 1\rho \log \M (-2  \rho \gl)
\, \le \, (1-\rho) \max_{\varrho\in [\rho,1]} \left\vert \frac{\dd}{\dd\varrho}g_\gl (\varrho)\right\vert
\end{equation}
where $g_\gl (\varrho):= \varrho ^{-1} \log \M ( -2\gl \varrho)$. Since
$\gl \in (0, \gl_0]$, one directly verifies that the rightmost term in 
\eqref{eq:nonGauss-2} is bounded by $C(\gl_0) \gl ^2$,
where $C(\gl_0)$ is a positive constant.

Another point in which the Gaussian character of the disorder enters 
is in the {\sl shifting procedure} of \eqref{eq:shift-1} and \eqref{eq:Hold} that gives \eqref{eq:expr-ub} and \eqref{eq:ai}. 
The shift has to be replaced by a {\sl tilt}
and one reduces to estimates that, except for constants that depend
on the law of the disorder, are the same as in the Gaussian case
(these steps are fully detailed for another model in \cite[Appendix~A.1 and Section 3]{cf:DGLT} and it is not very difficult to adapt them to this proof).
\qed

\section*{Acknowledgments}
The authors acknowledge the support of ANR, grant POLINTBIO.
T.B. and F.L.T. acknowledge also the support of  ANR,  grant  LHMSHE.

\end{document}